\documentclass[12pt, reqno]{amsart}
\usepackage{amsmath, amsthm, amscd, amsfonts, amssymb, graphicx}
\usepackage[bookmarksnumbered, colorlinks, plainpages]{hyperref}
\input{mathrsfs.sty}

\newtheorem{theorem}{Theorem}[section]
\newtheorem{lemma}[theorem]{Lemma}
\newtheorem{proposition}[theorem]{Proposition}
\newtheorem{corollary}[theorem]{Corollary}
\theoremstyle{definition}
\newtheorem{definition}[theorem]{Definition}
\newtheorem{example}[theorem]{Example}

\theoremstyle{remark}
\newtheorem{remark}[theorem]{Remark}
\numberwithin{equation}{section}
\newcommand{\ip}[2]{\left\langle#1,#2\right\rangle}
\newcommand{\iip}[2]{\left(#1|#2\right)}
\newcommand{\la}{\langle}
\newcommand{\ra}{\rangle}
\newcommand{\A}{{\mathscr A}}
\newcommand{\X}{{\mathscr X}}

\begin{document}
\title[Cauchy--Schwarz inequality]{A treatment of the Cauchy--Schwarz inequality in $C^*$-modules}

\author[Lj. Aramba\v si\' c, D. Baki\' c, M.S. Moslehian]{Ljiljana Aramba\v si\' c $^1$, Damir Baki\' c $^2$ and Mohammad Sal Moslehian $^3$}

\address{$^{1}$ $^{2}$ Department of Mathematics, University of Zagreb, Bijeni\v cka cesta 30, 10000 Zagreb, Croatia.}
\email{arambas@math.hr}
\email{bakic@math.hr}

\address{$^3$ Department of Pure Mathematics, Center of Excellence in
Analysis on Algebraic Structures (CEAAS), Ferdowsi University of Mashhad, P.O. Box
1159, Mashhad 91775, Iran.} \email{moslehian@ferdowsi.um.ac.ir and
moslehian@ams.org}



\begin{abstract} We study the Cauchy--Schwarz and some related inequalities in a semi-inner product module over a $C^*$-algebra ${\mathscr A}$.
The key idea is to consider a semi-inner product ${\mathscr A}$-module as a semi-inner product ${\mathscr A}$-module with respect to another semi-inner product. In this way, we improve some inequalities such as the Ostrowski inequality and an inequality related to the Gram matrix. The induced semi-inner products are also related to the the notion of covariance and variance. Furthermore, we obtain a sequence of nested inequalities that emerges from the Cauchy--Schwarz inequality. As a consequence, we derive some interesting operator-theoretical corollaries. In particular, we show that the sequence arising from our construction, when applied to a positive invertible element of a $C^*$-algebra, converges to its inverse.
\end{abstract}

\maketitle


\section{Introduction}

Let  $\A$ be a $C^*$-algebra. A \emph{(right) semi-inner product $\A$-module} is a
linear space ${\mathscr X}$ which is a right $\A$-module
with a compatible scalar multiplication ($\lambda(xa)=x(\lambda
a)=(\lambda x)a$ for all $x \in {\mathscr X},a \in \A,
\lambda \in {\mathbb C}$) endowed with an $\A$-semi-inner
product $\ip{\cdot}{\cdot}:{\mathscr X} \times{\mathscr X}\to
\A$ such that for all $x,y,z\in{\mathscr
X},\lambda\in{\mathbb C},a\in\A$, it holds
\begin{enumerate}
    \item[(i)] $\langle x,x\rangle\geq 0$;
    \item[(ii)] $\langle x, \lambda y + z\rangle=\lambda \langle x,y\rangle+
\langle x,z\rangle$;
    \item[(iii)] $\langle x,ya\rangle=\langle x,y\rangle a$;
    \item[(iv)] $\langle x,y\rangle^*=\langle y,x\rangle$.
\end{enumerate}
Obviously, every semi-inner product space is a semi-inner product
${\mathbb C}$-module. We can define a semi-norm on ${\mathscr X}$ by
$\|x\|=\|\la x,x\ra\|^\frac{1}{2}$, where the latter norm denotes
that in the $C^*$-algebra $\A$. A \emph{pre-Hilbert
$\A$-module} (or an \emph{inner-product module}) is a
semi-inner product module over $\A$ in which $\|\cdot\|$
defined as above is a norm. A pre-Hilbert $\A$-module $X$ such that $(X,\|\cdot\|)$ is complete is called a \emph{Hilbert $C^*$-module}. Each $C^*$-algebra $\A$
can be regarded as a Hilbert $\A$-module via $\langle a, b\rangle = a^*b\,\,(a, b \in \A)$.
Throughout the paper, $\widetilde{\A}$ stands for the
minimal unitization of $\A$. By $e$ we denote the unit
in $\widetilde{\A}$. If ${\mathscr X}$ is an $\A$-module then it can be regarded as an $\widetilde\A$-module
via $xe=x$.
The basic theory of $C^*$-algebras and Hilbert $C^*$-modules can be found in \cite{LAN, TAK, WO}.

One of the fundamental inequalities in a semi-inner product module $(\X,\ip{\cdot}{\cdot})$ over a $C^*$-algebra $\A$ is the Cauchy--Schwarz inequality. It states that
\begin{equation}\label{c-s}
\langle x,y\rangle \langle y,x\rangle  \leq \|y\|^2\,\langle
x,x\rangle\qquad(x,y\in \X)
\end{equation}
which generalizes the classical Cauchy--Schwarz inequality. Today many generalizations of the classical Cauchy--Schwarz inequality for integrals, isotone functionals as well as in the setting of inner product spaces are well-studied; see the book \cite{DRA}. Moreover, Niculescu \cite{NIC} and Joi\c{t}a \cite{JOI} have investigated the reverse of the Cauchy--Schwarz inequalities in the framework of $C^*$-algebras and Hilbert $C^*$-modules, see also \cite{M-P} and references therein. We also refer to another interesting paper by Ili\v sevi\'c and Varo\v{s}anec \cite{Di-Sa} of this type. Some operator versions of the Cauchy--Schwarz inequality with simple conditions for the case of equality are presented in \cite {FUJ}.

Using the Cauchy--Schwarz inequality in a new suitably defined semi-inner product on ${\mathscr X}$ we improve some known inequalities in semi-inner product modules. The most interesting improvement is the one of the Cauchy--Schwarz inequality itself. In this way, we improve some inequalities such as the Ostrowski inequality (see \cite{A-R}) and show that the Gram matrix $\left[\ip{x_i}{x_j}\right]$, where $x_1,\ldots, x_n\in\X$, is greater then or equal to some positive element of $M_n({\mathscr A})$. The induced semi-inner products are also related to the the notion of covariance and variance.

In the last section of the paper we repeat this technique by starting with an induced semi-inner product. This leads to a (possibly finite) sequence of nested inequalities refining the Cauchy--Schwarz inequality. Namely, for every $z\in{\mathscr X}$
such that $a:=\ip{z}{z}\neq 0$, we construct an increasing sequence of positive elements $(p_m(a))_m\in\A$ such that
$\langle x,x\rangle\ge \ip{x}{z} p_m(a)\ip{z}{x}$ for all $x,y\in {\mathscr X}$ (see Theorem~\ref{induct});
thereby, for $m=0$ we get $\|z\|^2\langle x,x\rangle\ge \ip{x}{z}\ip{z}{x}$, i.e. the Cauchy--Schwarz inequality.
By analyzing the sequence $(p_m(a))_m$ we obtain some interesting operator-theoretical consequences.
In particular, if $a$ is invertible, then $(p_m(a))_m$ converges in norm to $a^{-1}$ (Theorems~\ref{inv} and \ref{invinfty}).
Moreover, in Proposition~\ref{closed} we show that, for a positive operator $a$ on a Hilbert space ${\mathscr H}$, the sequence $(ap_m(a))_m$ converges in norm to the orthogonal projection to $\overline{\mbox{Im}\,a}$ if and only if $\mbox{Im}\,a$
is a closed subspace of ${\mathscr H}.$

\section{An induced family of semi-inner products}  \label{main-res}

Let $(\X,\ip{\cdot}{\cdot})$ be a semi-inner product module over a $C^*$-algebra $\A$.
For an arbitrary $z\in\X$ we define
\begin{equation}\label{z-inner}
\la\cdot,\cdot\ra_z:{\mathscr X}\times{\mathscr X}\to\A,\quad \la x,y\ra_z:= \|z\|^2\,\langle x,y\rangle-\langle x,z\rangle \langle z,y\rangle.
\end{equation}
It is easy to see that $\la \cdot,\cdot\ra_z$ is a semi-inner product on $\X.$ (Note that the case when $\la z,z\ra=0$ gives a trivial semi-inner product; however, this does not contradict the definition of the semi-inner product.)

In this section we show how one can improve, by using this new class of induced semi-inner products, several results known from the literature that rely on the Cauchy--Schwarz inequality for the original semi-inner product in an arbitrary module.


\subsection{The Gram matrix}
We begin with the Gram matrix $\left[\ip{x_i}{x_j}\right]$ in the matrix $C^*$-algebra $M_n(\A)$ of all $n \times n$ matrices
with entries from $\A$ due to it naturally appears in this context. Namely, positivity of the Gram matrix for two elements is strongly related to the Cauchy--Schwarz inequality.
To see this, let us write the Cauchy--Schwarz inequality in a matrix form.
Recall that a matrix $\left[\begin{array}{cc}a&b\\b^*&
c\end{array}\right]\in M_2(\A)$ with invertible $c\in
\A$ (resp. $a\in\A$) is positive if and only if
$a \geq 0, c \geq 0$ and $bc^{-1}b^* \leq a$ (resp. $a \geq 0, c
\geq 0$ and $b^*a^{-1}b \leq c$); see \cite{CHO}. Therefore,  \eqref{c-s} can be written as
\begin{eqnarray}\label{mat_intr1}
\left[\begin{array}{cc}\la x,x\ra & \la x,y\ra\\
\la x,y\ra^*& \|\la y,y\ra\|e\end{array}\right]
\geq 0\,,
\end{eqnarray}
where $e\in\widetilde{\A}$ is the unit. Since
$$
\left[\begin{array}{cc}\la x,x\ra & \la x,y\ra\\
\la x,y\ra^*& \|\la y,y\ra\|e\end{array}\right]
\geq
\left[\begin{array}{cc}\la x,x\ra & \la x,y\ra\\
\la x,y\ra^*& \la y,y\ra\end{array}\right]
\geq 0,$$
it follows that positivity of the Gram matrix sharpens the Cauchy--Schwarz inequality.

\begin{remark}
A number of arguments can be simplified if we use positivity of the Gram matrix. For example, it was proved in \cite[Theorem~2.1]{Di-Sa} that for $x,y\in\X$ such that $|y|$ belongs to the center of $\A,$ a stronger version of the Cauchy--Schwarz inequality holds, namely, $\la x,y\ra\la y,x\ra\le \la x,x\ra\la y,y\ra.$
From positivity of the Gram matrix it follows that for every $x,y\in\X$ and every $\varepsilon>0$ we have
$$\left[\begin{array}{cc}\la x,x\ra & \la x,y\ra\\
\la x,y\ra^*& \la y,y\ra+\varepsilon e\end{array}\right]
\geq 0,$$
or, equivalently, $\la x,y\ra(\la y,y\ra+\varepsilon e)^{-1}\la y,x\ra\le \la x,x\ra.$ If $|y|$ belongs to the center of $\A,$ we get, by multiplying by $(\la y,y\ra+\varepsilon e)^{\frac{1}{2}}$ on both sides, $\la x,y\ra\la y,x\ra\le \la x,x\ra(\la y,y\ra+\varepsilon e).$ Since $\varepsilon>0$ is arbitrary, we have
$$\la x,y\ra\la y,x\ra\le \la x,x\ra\la y,y\ra.$$

Because of this stronger version of the Cauchy--Schwarz inequality, we can define another family of semi-inner products on $\X$: for every $z\in\X$ such that $|z|$ belongs to the center of $\A,$ the mapping
$$\{\cdot,\cdot\}_z:{\mathscr X}\times{\mathscr X}\to\A,\quad \{x,y\}_z:= \la z,z\ra \la x,y\ra-\la x,z\ra \la z,y\ra$$
is a semi-inner product on $\X.$
\end{remark}

\vspace{0.1in}

Let us now consider the Gram matrix $\left[\ip{x_i}{x_j}\right]\in M_n(\A)$ for an arbitrary number of elements $x_1,\ldots, x_n$ in a semi-inner product module $(\X,\ip{\cdot}{\cdot})$.
It is known that $[\la x_i,x_j\ra]\ge 0,$ i.e., the Gram matrix is a positive element of the
$C^*$-algebra $M_n(\A)$. The interesting fact about this inequality is that it is self-improving, as we show in the following theorem.

\begin{theorem}\label{p1}
Let $\A$ be a $C^*$-algebra and $({\mathscr
X},\ip{\cdot}{\cdot})$  a semi-inner product $\A$-module. Let $n\in{\mathbb N}$ and $x_1, \ldots, x_n \in
{\mathscr X}$. Then for every $z\in {\mathscr X}$ we have
\begin{equation}\label{bas}
   \|z\|^2\big[\ip{x_i}{x_j}\big]\ge \big[\ip{x_i}{z}\ip{z}{x_j}\big].
 \end{equation}
\end{theorem}
\begin{proof} Let us first prove that $[\ip{x_i}{x_j}]$ is
positive in $M_n(\A)$ (the proof is included for the convenience of the
reader, see \cite[Lemma~4.2]{LAN}). Since
$\ip{\cdot}{\cdot}$ is a semi-inner product on ${\mathscr X}$ it
holds that
$$\ip{\sum_{i=1}^nx_ia_i}{\sum_{i=1}^nx_ia_i}\ge
0,\quad (a_1,\ldots,a_n\in\A).$$ Then
\begin{equation}\label{pos}
\sum_{i,j=1}^n a_i^*\ip{x_i}{x_j}a_j\ge
0,\quad(a_1,\ldots,a_n\in\A).
\end{equation}
By \cite[Lemma IV.3.2]{TAK} we know that a matrix $[c_{ij}]\in
M_n(\A)$ is positive if and only if
$\sum_{i,j=1}^na_i^*c_{ij}a_j \geq 0$ for all $a_1, \ldots, a_n \in
\A$. Therefore, \eqref{pos} means that the matrix
$[\ip{x_i}{x_j}]$ is positive.

It holds for an arbitrary semi-inner product, so, choosing  $\ip{\cdot}{\cdot}_z$ instead of $\ip{\cdot}{\cdot}$, we get $\big[\la x_i,x_j\ra_z\big]\ge 0$, which is exactly \eqref{bas}.
\end{proof}

The following corollary is a direct consequence of the preceding theorem. A positive linear mapping $\Phi:\A \to {\mathscr B}$, where ${\mathscr B}$ is a $C^*$-subalgebra of $\A$, is called a \textit{left multiplier} if $\Phi(ab)=\Phi(a)b\,\,(a \in \A, b \in {\mathscr B}$).

\begin{corollary}\label{pos-mult} Let $(\X,\la\cdot,\cdot\ra)$ be a semi-inner product $\A$-module, ${\mathscr B}$ a $C^*$-subalgebra of $\A$ and $\Phi: \A \to {\mathscr B}$ a positive left multiplier. Then
\begin{eqnarray}\label{cond}
\|\Phi(\la z,z\ra)\|[\Phi(\la x_i,x_j\ra)] \geq [\Phi(\la x_i,z\ra )\Phi(\la z,x_j)]
\end{eqnarray}
for all $x_1,\cdots,x_n,z\in\X$.
\end{corollary}
\begin{proof}
Given a left multiplier $\Phi:\A \to {\mathscr B}$, any semi-inner product $\A$-module $\X$ becomes a semi-inner product ${\mathscr B}$-module with respect to
\begin{eqnarray}\label{multip}
[x,y]_\Phi=\Phi(\langle x, y\rangle),\quad (x,y\in{\mathscr X}).
\end{eqnarray}
By \eqref{bas}, it holds
   $$\|[z,z]_\Phi\|\big[[x_i,x_j]_\Phi\big]\ge \big[ [x_i,z]_\Phi [z,x_j]_\Phi\big]\,.$$
\end{proof}
\vspace{.1in}

\begin{remark}
Let ${\mathscr X}$ be a $C^*$-algebra regarded as
a Hilbert $C^*$-module over itself. Since every conditional
expectation $\Phi: \A \to {\mathscr B}$ is a completely positive left multiplier
(cf. \cite[IV, \S 3]{TAK}), the preceding corollary is an extension of \cite[Theorem 1]{B-D} for conditional
expectations.
\end{remark}

\subsection{A covariance--variance inequality}

Another application of Theorem~\ref{p1} is the
covariance--variance inequality in semi-inner product
$C^*$-modules. The interested reader is referred
to \cite{B-D,F-F-N-T, LIN} for some generalizations of
covariance--variance inequality. Let us begin with a definition and some known
examples.

\begin{definition} \label{cov} Let $\A$ be a $C^*$-algebra, $({\mathscr X},\ip{\cdot}{\cdot})$ be a semi-inner
product $\A$-module and $x,y,z\in {\mathscr X}$. The
\emph{covariance $\textrm{cov}_z(x,y)$ between $x$ and $y$ with
respect to $z$} is defined to be the element $\la x,y\ra_z$ of
$\A$. The element $\textrm{cov}_z(x,x)$ is said to be the
\emph{variance of $x$ with respect to $z$} and denoted by
$\textrm{var}_z(x)$.
\end{definition}
\begin{example} Given a Hilbert space ${\mathscr H}$, vectors $x,y\in{\mathscr
H}$ and operators $S, T \in {\mathbb B}({\mathscr H})$, covariance
and variance of operators was defined in \cite{LIN} as
$$\textrm{cov}_{x,y}(S,T)=\|y\|^2(Sx|Tx)-(Sx|y)(y|Tx).$$ Observe that $\textrm{cov}_{x,y}(S,T)=\iip{Sx}{Tx}_y$. In
the case where $\|x\|=1$ and $y=x$ we get the notion of covariance
of two operators $T$ and $S$ introduced in \cite{F-F-N-T} as
$$\textrm{cov}_x(S,T)=(Sx|Tx)-(Sx|x)(x|Tx)\,.$$
A notion of covariance and variance of Hilbert space operators was investigated in
\cite{F-F-N-T, SEO}. In addition, Enomoto
\cite{ENO} showed a close relation of the operator
covariance--variance inequality with the Heisenberg uncertainty
principle and pointed out that it is exactly the generalized
Schr\"odinger inequality. For more information on related ideas and concepts we refer the reader
to \cite[Section 5]{PAR}.

Another remarkable fact is that for a unit vector $x\in {\mathscr H}$, the determinant of the positive semidefinite Gram matrix
\[\left[\begin{array}{ccc} \iip{Sx}{Sx}& \iip{Sx}{Tx}& \iip{Sx}{x}\\
\iip{Tx}{Sx}& \iip{Tx}{Tx}& \iip{Tx}{x}\\
\iip{x}{Sx}& \iip{x}{Tx}& \iip{x}{x}\end{array}\right]\] is the
difference $\textrm{var}_x(S)\textrm{var}_x(T) -
|\textrm{cov}_x(S,T)|^2$ and is nonnegative; see \cite{F-I-N-S}.
\end{example}
\begin{example}
Recall that if $(\Omega,\mu)$ is a probability measure space, then
$Ef=\int_\Omega fd\mu$ is the expectation of the random variable
$f\in L^2(\Omega,\mu)$. Then the covariance between $f$ and $g$ is
defined to be $\textrm{cov}(f,g)=E(\overline{f}g)-\overline{Ef}\,Eg$
and variance of $f$ is $\textrm{cov}(f,f)$. We can obtain this by
considering $L^2(\Omega,\mu)$ as a Hilbert ${\mathbb C}$-module via
the usual inner product $\langle
f,g\rangle=\int_\Omega\overline{f}g$.
\end{example}


Let $\A$ be a $C^*$-algebra and ${\mathscr X}$ be a
semi-inner product $\A$-module.
The Cauchy--Schwarz inequality for $\textrm{cov}_z(\cdot,\cdot)$
is known as the \emph{covariance--variance inequality}.
Therefore, Theorem~\ref{p1} can be also stated in the following form.

\begin{theorem}[Generalized covariance-variance inequality]
\label{gv-ci} Let $\A$ be a $C^*$-algebra and ${\mathscr
X}$ be a semi-inner product $\A$-module. Let $x_1, \ldots,
x_n,z \in {\mathscr X}$. Then the matrix
$[\textrm{cov}_z(x_i,x_j)]\in M_n(\A)$ is positive.
\end{theorem}


\vspace{.1in}
Assume that $\A$ is a $C^*$-algebra acting on a Hilbert
space, ${\mathscr B}$ is one of its $C^*$-subalgebras and ${\mathscr X}$ is a
Hilbert $\A$-module. Let us fix  a positive left
multiplier mapping $\Phi$ and $x \in {\mathscr X}$ such that $\|\Phi(\langle x,x\rangle)\|=1$.  For operators $A$ and $B$ in the algebra ${\mathbb B}({\mathscr X})$ of all adjointable operators
on ${\mathscr X}$ we could define the covariance
of $A, B$ and variance of $A$ by
$$\textrm{cov}(A,B)=\Phi(\langle
Ax, Bx\rangle)-\Phi(\langle Ax, x\rangle)\Phi(\langle x,
Bx\rangle)$$ and $\textrm{var}(A)=\textrm{cov}(A,A)$,
respectively. Observe that, if we regard $X$ as a
semi-inner product $\mathscr{A}$-module with respect to
$[\cdot,\cdot]_\Phi$ defined by \eqref{multip}, then we have
$\textrm{cov}(A,B)=\textrm{cov}_x(Ax,Bx).$ Therefore,
\begin{eqnarray*}
\left[\begin{array}{cc}\textrm{var}(A)&\textrm{cov}(A,B)\\\textrm{cov}(A,B)^*&
\textrm{var}(B)\end{array}\right]=\left[\begin{array}{cc}\textrm{var}_x(Ax)&\textrm{cov}_x(Ax,Bx)\\\textrm{cov}_x(Ax,Bx)^*&
\textrm{var}_x(Bx)\end{array}\right]\ge 0.
\end{eqnarray*}


\subsection{An Ostrowski type inequality}

Here we show that some Ostrowski-type inequalities can be viewed as the Cauchy--Schwarz inequality with respect to a new semi-inner product.

It was proved in \cite{Ma} that for any three elements in a real inner product space $(H,(\cdot|\cdot))$ it holds
$$\left\vert\|z\|^2( x|y)-(x|z)(y|z)\right\vert^2
\leq \left( \|z\|^2\|x\|^2-(x|z)^2\right)\left(\|z\|^2\|y\|^2-(y|z)^2\right).$$
Since here $H$ is a real vector space, this may be written as
$$|(x|y)_z|^2\leq (x|x)_z (y|y)_z$$
and this is exactly the Cauchy--Schwarz inequality for $(\cdot|\cdot)_z.$
Therefore, the Cauchy--Schwarz inequality for a semi-inner product $\la \cdot,\cdot\ra_z$ on a semi-inner product module $\X$, i.e.
\begin{eqnarray}\label{ma}
 \nonumber
(\|z\|^2\,\langle y,x\rangle-\langle y,z\rangle \langle z,x\rangle)(\|z\|^2\,\langle x,y\rangle-\langle x,z\rangle \langle z,y\rangle)\\
\leq \left\Vert\|z\|^2\,\langle x,x\rangle-\langle x,z\rangle \langle z,x\rangle\right\Vert (\|z\|^2\,\langle y,y\rangle-\langle y,z\rangle \langle z,y\rangle)
\end{eqnarray}
generalizes the result from \cite{Ma}. In the special case when $\la x,z\ra=0$ we get
\begin{eqnarray}\label{LjR}
|\ip{z}{y}|^2\le \frac{\|z\|^2}{\|x\|^2}(\|x\|^2\vert
y\vert^2-|\ip{x}{y}|^2),
\end{eqnarray}
which is the Ostrowski inequality in a semi-inner product $C^*$-module (see \cite{A-R}).
Since \eqref{ma} is the Cauchy--Schwarz inequality and \eqref{LjR} is its special case, Theorem~\ref{p1} improves both of them. Namely,
$$\left[\begin{array}{cc}\la x,x\ra_z & \la x,y\ra_z\\
\la x,y\ra_z^*& \la y,y\ra_z\end{array}\right]
\geq 0,$$
(which is exactly \eqref{bas} for $n=2$) improves \eqref{ma}, and it improves \eqref{LjR} in the case $\la x,z\ra=0$.


\section{A nested sequence of inequalities}

In this section we show that the Cauchy--Schwarz inequality can be improved by a sequence of nested
inequalities. To do that, let us first fix some notation.

Let $\A$ be a $C^*$-algebra and $e\in\widetilde\A$. For a positive element $a \in \A$, $a\not =0$, define
\begin{equation}\label{efovi}
\begin{array}{ll}
f_0(a)=a,&g_0(a)=\|f_0(a)\|e-f_0(a),\\
f_1(a)=f_0(a)g_0(a),&g_1(a)=\|f_1(a)\|e-f_1(a),\\
\ldots& \\
f_m(a)=f_{m-1}(a)g_{m-1}(a),&g_m(a)=\|f_m(a)\|e-f_m(a),\\
\ldots
\end{array}
\end{equation}

Observe that all $f_m(a)$ and $g_m(a)$ are polynomials in $a$. An
easy inductive argument shows that all $f_m(a)$ and $g_m(a)$ are
positive elements as well, and for all $m\ge 0$  it holds
\begin{equation}\label{svojstva}
f_{m+1}(a)=f_m(f_1(a)),\quad g_{m+1}(a)=g_m(f_1(a)).
\end{equation}
It may happen that $f_m(a)=0$ for some
$m\in \Bbb N$ (see Proposition \ref{kon} below); then, obviously,
$f_k(a)=0$ for all $k \geq m$. On the other hand, if $f_m(a) \not
=0$ for some $m$, then, by definition, $f_j(a) \not =0,\,\forall j
\leq m$. Thus, for each $m$ such that $f_m(a) \not =0$ we can define
\begin{equation}\label{peovi}
\begin{array}{l}
p_0(a)=\frac{e}{\|f_0(a)\|},\\
p_1(a)=\frac{e}{\|f_0(a)\|}+\frac{g_0(a)^2}{\|f_0(a)\|\cdot \|f_1(a)\|},\\
p_2(a)=\frac{e}{\|f_0(a)\|}+\frac{g_0(a)^2}{\|f_0(a)\|\cdot \|f_1(a)\|}+\frac{g_0(a)^2g_1(a)^2}{\|f_0(a)\|\cdot \|f_1(a)\|\cdot \|f_2(a)\|},\\
\ldots \\
p_m(a)=\frac{e}{\|f_0(a)\|}+
\sum_{l=1}^m\left(\frac{1}{\prod_{k=0}^l\|f_k(a)\|}\prod_{k=0}^{l-1}g_k(a)^2\right).
\end{array}
\end{equation}
It is convenient here to make the following convention: if $m$ is the last index such that $f_m(a) \not =0$ then we define
\begin{equation}\label{peovi2}
p_j(a)=p_m(a),\quad (j> m).
\end{equation}
Thus, we can treat $(p_m(a))$ as an infinite sequence of
positive elements in $\A$ even in the case when there is $m\ge 0$ such that $f_m(a)=0.$

\begin{remark}
It is obvious that $0 \leq p_0(a) \leq p_1(a) \leq \ldots \leq p_m(a) \leq \ldots$.
Observe also that $p_m(a)$'s which are defined by \eqref{peovi} are all different. Indeed, suppose $p_{m-1}(a)$ and $p_{m}(a)$ are defined by \eqref{peovi} and $p_{m-1}(a)=p_{m}(a).$ Then $$\frac{1}{\prod_{k=0}^m\|f_k(a)\|}\prod_{k=0}^{m-1}g_k(a)^2=0$$ which implies  $g_0(a)g_1(a)\cdots g_{m-1}(a)=0$ and therefore
\begin{eqnarray*}
  f_m(a) &=& f_{m-1}(a)g_{m-1}(a)=f_{m-2}(a)g_{m-2}(a)g_{m-1}(a) \\
  &=& \ldots=f_0(a)g_0(a)\cdots g_{m-1}(a)=0.
\end{eqnarray*}
This is the contradiction, since $p_m(a)$ is defined by \eqref{peovi}.
\end{remark}

\begin{theorem} \label{induct}
Let ${\mathscr X}$ be a module over a $C^*$-algebra $\A $, and let $\ip{\cdot}{\cdot}$ be any $\A $-valued semi-inner product on ${\mathscr X}$.
For each $z\in{\mathscr X}$
such that $\ip{z}{z}\neq 0$ it holds
\begin{eqnarray*}\label{nest}
\langle x,x\rangle&\ge& \ldots\ge \ip{x}{z} p_m(\ip{z}{z})
\ip{z}{x}\ge \ip{x}{z} p_{m-1}(\ip{z}{z})\ip{z}{x}\\
&\ge&\ldots\ge \ip{x}{z} p_0(\ip{z}{z})
\ip{z}{x}=\frac{1}{\|z\|^2}\ip{x}{z}\ip{z}{x}\ge 0.
\end{eqnarray*}
\end{theorem}

\begin{proof}
We will prove by induction that
\begin{equation}\label{jedan}
\la x,x\ra_*\ge \la x,z\ra_* p_m(\la
z,z\ra_*) \la z,x\ra_*
\end{equation}
holds true for all $m\geq 0$, for each $z \in {\mathscr X}$ and for every $\A $-valued semi-inner product
$\ip{\cdot}{\cdot}_*$ on ${\mathscr X}$ such that $\ip{z}{z}_* \not =0$.

For $m=0$ this is precisely the statement of Theorem \ref{p1} for $n=1$.

Suppose that (\ref{jedan}) holds for some $m$ and for all $z$ and
$\ip{\cdot}{\cdot}_*$ such that $\ip{z}{z}_* \not =0$. Choose an arbitrary
semi-inner product $\ip{\cdot}{\cdot}$ on ${\mathscr X}$ such that  $\ip{z}{z} \not =0$.
If $f_{m+1}(\ip{z}{z})=0$, there is nothing to prove since then, by our convention, $p_{m+1}(\ip{z}{z})=p_{m}(\ip{z}{z})$.

Suppose now that $f_{m+1}(\la z,z\ra)\neq 0.$ Then, by \eqref{svojstva}, $f_{m}(f_1(\la z,z\ra)=f_{m}(\la z,z\ra_z)\neq 0.$
By the inductive assumption (for the semi-inner product $\la
\cdot,\cdot\ra_z$) it holds
$$\la x,x\ra_z\ge \la x,z\ra_z p_m(\la
z,z\ra_z) \la z,x\ra_z,$$
that is,
\begin{equation}\label{mm}
\|z\|^2\ip{x}{x}\ge \ip{x}{z}\ip{z}{x}+ \la
x,z\ra_z p_m(\la z,z\ra_z) \la z,x\ra_z.
\end{equation}
Observe that $\la z,z\ra_z=f_1(\ip{z}{z})$, so $\|z\|^2e-\ip{z}{z}$
and $p_m(\la z,z\ra_z)$ commute. Therefore
\begin{eqnarray*}
  \la x,z\ra_z p_m(\la z,z\ra_z) \la z,x\ra_z &=& (\|z\|^2\ip{x}{z}-\ip{x}{z}\ip{z}{z}) p_m(\la z,z\ra_z) \\
 &  &\cdot(\|z\|^2\ip{z}{x}-\ip{z}{z}\ip{z}{x})\\
 &=& \ip{x}{z}(\|z\|^2e-\ip{z}{z})^2p_m(\la z,z\ra_z)\ip{z}{x}\\
 &=& \ip{x}{z}\left(g_0(\la z,z\ra)\right)^2p_m\left(f_1(\la z,z\ra)\right)\ip{z}{x}.
\end{eqnarray*}
Since $f_{m+1}(\la z,z\ra)\neq 0$ and $f_{m}(\la z,z\ra_z)\neq 0,$ the elements $p_{m+1}(\la z,z\ra)$ and $p_{m}(\la z,z\ra_z)$ are defined by \eqref{peovi}. It is easy to verify, using \eqref{efovi} and \eqref{svojstva}, that
\begin{equation}\label{a_p}
e+g_0(\la z,z\ra)^2p_m(f_1(\la z,z\ra))=\|z\|^2p_{m+1}(\ip{z}{z}),
\end{equation}
which, together with \eqref{mm}, gives $\ip{x}{x}\ge
\ip{x}{z}p_{m+1}(\ip{z}{z}) \ip{z}{x}$.

To complete the proof it only remains to recall that $p_{m}(\ip{z}{z})\ge p_{m-1}(\ip{z}{z})\ge
\ldots\ge p_{0}(\ip{z}{z})\ge 0$, for every $z\in{\mathscr X}$.
\end{proof}

If $\ip{z}{z}$ is not a scalar multiple of the unit, then $f_1(\ip{z}{z})\neq 0$ and the preceding theorem
strictly refines the inequality from Theorem
\ref{p1}. Moreover, if $f_m(\ip{z}{z})\neq 0$ for all $m \in \Bbb
N$, Theorem \ref{induct} provides an infinite sequence of
inequalities.
On the other hand, if $f_m(\ip{z}{z})=0$ for some $m\ge 0$, then, by \eqref{peovi2}, only finitely
many inequalities are obtained. The following
proposition characterizes all such elements $z \in {\mathscr X}$. It
turns out that the sequence of inequalities obtained in Theorem
\ref{induct} is finite precisely when $\ip{z}{z}$ has a finite
spectrum.

\begin{proposition}\label{kon}
Let $a$ be a positive  element of a $C^*$-algebra
$\A\subseteq\mathbf{B}({\mathscr H}).$ Then there exists $m\in\mathbb{N}$ such that $f_{m}(a)=0$
if and only if $a$ has a finite spectrum.
\end{proposition}
\begin{proof}
Suppose that there is
$m\in\mathbb{N}$ such that $f_{m}(a)=0$. Let $\lambda \in \sigma(a)$.
Then $f_{m}(\lambda)\in f_{m}(\sigma(a))=\sigma(f_{m}(a))=\{0\}$.
This shows that $\sigma(a)$ is contained in a finite set, namely in the set of all zeros of the polynomial $f_{m}$.

To prove the converse, suppose that $\sigma(a)$ is a finite set. First observe that $0\in\sigma(f_m(a))$ for all $m\ge 1.$
Let $m\ge 1$ be such that $f_m(a)\neq 0$ (if such $m$ does not exist, we are done). Since
$\sigma(f_{m+1}(a))=\{\|f_m(a)\|\lambda-\lambda^2:\lambda\in\sigma(f_m(a))\},$ and since 0 and $\|f_m(a)\|$ are two different elements of $\sigma(f_m(a))$ such that $\|f_m(a)\| 0-0^2=\|f_m(a)\|\cdot \|f_m(a)\|-\|f_m(a)\|^2=0$, we conclude that $\textup{card}\,\sigma(f_{m+1}(a))<\textup{card}\,\sigma(f_{m}(a)).$ Since $\sigma(a)$ is finite, there is $m$ such that $\sigma(f_m(a))=\{0\},$ i.e. $f_m(a)=0$.
\end{proof}

\vspace{.1in}

\begin{remark}\label{spaces}
Suppose that $\A=\Bbb C$, i.e.~that ${\mathscr X}$ is a
semi-inner product space. Then for each $z \in {\mathscr X}$ the
spectrum $\sigma(\ip{z}{z})$ is a singleton, so
$f_1(\ip{z}{z})=0$. Hence, in this situation, the sequence of
inequalities from Theorem \ref{induct} terminates already at the
first step. In other words, Theorem \ref{induct} reduces then to
Theorem \ref{p1}. Therefore, Theorem \ref{induct} gives us a new
(possibly finite) sequence of inequalities only if the underlying
$C^*$-algebra is different from the field of complex numbers.
\end{remark}

Let us first consider the case from the preceding proposition, when the sequence of the inequalities is finite. The following result is interesting in its own. If $a \in \A$ is such that $f_M(a)\not =0$ and $f_{M+1}(a)=0$ for some $M \in \Bbb N$, we show that,
roughly speaking, $p_M(a)$ is the inverse of $a$.

\begin{theorem}\label{inv}
Let $a\neq 0$ be a positive element in a $C^*$-algebra $\A \subseteq \Bbb B({\mathscr H})$ with a finite spectrum. Let $M$ be the number with the property $f_M(a)\not =0$, $f_{M+1}(a)=0$.
Then $ap_M(a)$ is the orthogonal projection to the image of $a$. In particular, if $a$ is an invertible operator, then $p_M(a)=a^{-1}$.
\end{theorem}

\begin{proof}
Let us first observe that, since the spectrum of $a$ is finite,  $\text{Im}\,a$ is a closed subspace of ${\mathscr H}$.
For every $\lambda\in\mathbb{R}$ and $l\in\mathbb{N}$ it holds
\begin{eqnarray*}
  \lambda\prod_{k=0}^{l-1}g_k(\lambda)&=&f_0(\lambda)g_0(\lambda)g_1(\lambda)\cdots g_{l-1}(\lambda)=f_1(\lambda)g_1(\lambda)\cdots g_{l-1}(\lambda) \\
    &=&f_2(\lambda)g_2(\lambda)\cdots g_{l-1}(\lambda)=\ldots=f_{l}(\lambda).
\end{eqnarray*}
Since $f_M(a)\neq 0$ and $f_{M+1}(a)=0$ we conclude that $p_0(a),\ldots,p_M(a)$ are defined by \eqref{peovi}, while, by \eqref{peovi2}, $p_j(a)=p_M(a)$ for $j\ge M+1.$
Therefore, for $\lambda\neq 0$ and $m=1,\ldots,M$ it holds
\begin{equation}\label{pmfm}
p_m(\lambda)=\frac{1}{\|f_0(a)\|}+
\frac{1}{\lambda^2}\sum_{l=1}^m\frac{f_l(\lambda)^2}{\prod_{k=0}^l\|f_k(a)\|}.
\end{equation}

Let $m\in\{1,\ldots,M\}.$ Since $f_m(a)\geq 0$ (and $f_m(a)\neq 0$), $\|f_m(a)\|$ is the maximum of the set $\sigma(f_m(a))=f_m(\sigma(a))$. Let $\lambda_m\in\sigma(a)$ be such that $f_m(\lambda_m)=\|f_m(a)\|.$ Then $g_m(\lambda_m)=\|f_m(a)\|-f_m(\lambda_m)=0$ and therefore $f_j(\lambda_m)=0$ for all $j\ge m+1.$ Since obviously $\lambda_m \neq 0$, \eqref{pmfm} gives  $p_j(\lambda_m)=p_m(\lambda_m)$ for all $j\in \{m,\ldots,M\}$.
Therefore, for $j\in \{m,\ldots,M\}$ we have
$$p_j(\lambda_m)=p_m(\lambda_m)=\frac{1}{\|f_0(a)\|}+
\frac{1}{\lambda_m^2}\left(\sum_{l=1}^{m-1}\frac{f_l(\lambda_m)^2}{\prod_{k=0}^l\|f_k(a)\|}+\frac{f_m(\lambda_m)^2}{\prod_{k=0}^m\|f_k(a)\|}\right).$$
Using $f_m(\lambda_m)=\|f_m(a)\|$,  for all  $j\in \{m,\ldots,M\}$ we get
\begin{eqnarray}\label{for1}
p_j(\lambda_m)&=&\frac{1}{\|f_0(a)\|}+
\frac{1}{\lambda_m^2}\left(\sum_{l=1}^{m-1}\frac{f_l(\lambda_m)^2}{\prod_{k=0}^l\|f_k(a)\|}+\frac{f_m(\lambda_m)}{\prod_{k=0}^{m-1}\|f_k(a)\|}\right)\\
\nonumber &=&\frac{1}{\|f_0(a)\|}+
\frac{1}{\lambda_m^2}\left(\sum_{l=1}^{m-2}\frac{f_l(\lambda_m)^2}{\prod_{k=0}^l\|f_k(a)\|}+
\frac{f_{m-1}(\lambda_m)^2}{\prod_{k=0}^{m-1}\|f_k(a)\|}+\frac{f_m(\lambda_m)}{\prod_{k=0}^{m-1}\|f_k(a)\|}\right)\\
\nonumber &=&\frac{1}{\|f_0(a)\|}+
\frac{1}{\lambda_m^2}\left(\sum_{l=1}^{m-2}\frac{f_l(\lambda_m)^2}{\prod_{k=0}^l\|f_k(a)\|}+
\frac{f_{m-1}(\lambda_m)^2+f_m(\lambda_m)}{\prod_{k=0}^{m-1}\|f_k(a)\|}\right).
\end{eqnarray}
Observe that for every $k\in\mathbb{N}$ and every $\lambda\in\mathbb{R}$ it holds: $f_{k-1}(\lambda)^2+f_k(\lambda)=f_{k-1}(\lambda)^2+f_{k-1}(\lambda)g_{k-1}(\lambda)=
f_{k-1}(\lambda)(f_{k-1}(\lambda)+g_{k-1}(\lambda))=f_{k-1}(\lambda)\|f_{k-1}(a)\|$. Therefore, for $j\in \{m,\ldots,M\}$,
\begin{eqnarray}
\nonumber p_j(\lambda_m)&=&\frac{1}{\|f_0(a)\|}+
\frac{1}{\lambda_m^2}\left(\sum_{l=1}^{m-2}\frac{f_l(\lambda_m)^2}{\prod_{k=0}^l\|f_k(a)\|}+
\frac{f_{m-1}(\lambda_m)\|f_{m-1}(a)\|}{\prod_{k=0}^{m-1}\|f_k(a)\|}\right)\\
\label{for2}
&=&\frac{1}{\|f_0(a)\|}+
\frac{1}{\lambda_m^2}\left(\sum_{l=1}^{m-2}\frac{f_l(\lambda_m)^2}{\prod_{k=0}^l\|f_k(a)\|}+
\frac{f_{m-1}(\lambda_m)}{\prod_{k=0}^{m-2}\|f_k(a)\|}\right).
\end{eqnarray}
We now proceed recursively in the same way as (\ref{for2}) is obtained from (\ref{for1}) to get
\begin{eqnarray*}
p_j(\lambda_m)&=&\frac{1}{\|f_0(a)\|}+
\frac{1}{\lambda_m^2}\left(\sum_{l=1}^{1}\frac{f_l(\lambda_m)^2}{\prod_{k=0}^l\|f_k(a)\|}+
\frac{f_{2}(\lambda_m)}{\prod_{k=0}^{1}\|f_k(a)\|}\right)\\
&=&\frac{1}{\|f_0(a)\|}+
\frac{1}{\lambda_m^2}\left(\frac{f_1(\lambda_m)^2}{\|f_0(a)\|\|f_1(a)\|}+
\frac{f_{2}(\lambda_m)}{\|f_0(a)\|\|f_1(a)\|}\right)\\
&=&\frac{1}{\|f_0(a)\|}+
\frac{1}{\lambda_m^2}\left(\frac{f_1(\lambda_m)\|f_{1}(a)\|}{\|f_0(a)\|\|f_1(a)\|}\right)\\
&=&\frac{1}{\|f_0(a)\|}+\frac{1}{\lambda_m^2}\frac{\lambda_m g_0(\lambda_m)}{\|f_0(a)\|}\\
&=&\frac{\lambda_m+g_0(\lambda_m)}{\lambda_m\|f_0(a)\|}=\frac{\|f_0(a)\|}{\lambda_m\|f_0(a)\|}=\frac{1}{\lambda_m}
\end{eqnarray*}
for $j\in \{m,\ldots,M\}$.

After all, we have proved:
if $\lambda_m\in\sigma(a)$ is such that $\|f_m(a)\|=f_m(\lambda_m)$ for some $m\in\{0,\ldots,M\}$ then $p_j(\lambda_m)=\frac{1}{\lambda_m}$ for all $j\in \{m,\ldots,M\}$.

Let us take particular $\lambda\in\sigma(a),\lambda\neq 0.$ From $f_{M+1}(a)=0$ it follows that $f_{M+1}(\lambda)=0$. Then there exists $m\le M$  such that $f_m(\lambda)\neq 0$ and $f_{m+1}(\lambda)=0.$ Then from $f_{m+1}(\lambda)=f_m(\lambda)g_m(\lambda)$ we get $g_m(\lambda)=0,$ i.e., $f_m(\lambda)=\|f_m(a)\|.$ This means that for every $\lambda\in\sigma(a)$ there is $m\le M$ such that
$p_j(\lambda)=\frac{1}{\lambda}$ for all $j\in\{m,\ldots,M\}.$
Since $\sigma(a)$ is finite, there is $m\le M$ such that
$$p_{m}(\lambda)=\frac{1}{\lambda},\quad \forall \lambda\in\sigma(a)\setminus\{0\}$$
and therefore $$p_{M}(\lambda)=\frac{1}{\lambda},\quad \forall \lambda\in\sigma(a)\setminus\{0\}.$$
Then
$$\lambda p_{M}(\lambda)=
\left\{
\begin{array}{ll}
   1, & \lambda\in\sigma(a)\setminus\{0\}, \\
   0, & \lambda\in\sigma(a)\cap\{0\}.
\end{array}
\right.$$
This is precisely what we need to conclude that $a p_{M}(a)$ is the orthogonal projection to $\textup{Im}\, a.$ In the case when $a$ is invertible, then $\lambda p_{M}(\lambda)=1$ for all $\lambda\in\sigma(a)$, so $a p_{M}(a)=e.$
\end{proof}

\vspace{.1in}
Let us now consider the sequence $(p_m(a))$ in full generality. Again, we assume that $a$ is a positive operator on some Hilbert space ${\mathscr H}$ (i.e.~ $\A$ is represented faithfully on
${\mathscr H}$). Denote by $p$ the orthogonal projection to $\overline{\mbox{Im}\,a}$. If $a$ has a finite spectrum we have seen in the preceding theorem that $ap_M(a)$ is equal to $p$, where $M$ is less than or equal to the number of elements of $\sigma(a)$; in particular, if $a$ is an invertible operator, then $p_M(a)=a^{-1}$.

If $\sigma(a)$ is an infinite set we know that $f_m(a)$ is never equal to $0$; thus, $(p_m(a))$ is in this case an increasing sequence of positive elements of $\Bbb B({\mathscr H})$.
It would be natural to expect that in this situation the sequence $(ap_m(a))$ converges to $p$ in norm. However, this is not true in general, as the following example shows.

\begin{example}\label{compact}
Suppose that $a$ is a positive compact operator with an infinite spectrum. Then the sequence $(ap_m(a))$ cannot converge to $p$ in norm.

Indeed, suppose the opposite. Observe that $p_m(a) \in C^*(a)$, for all $m \geq 0$, where $C^*(a)$ denotes the $C^*$-algebra generated by $a$. Since $C^*(a)$ is closed, the assumption would imply $p\in C^*(a)$. But this is impossible: since $\sigma(a)$ is an infinite set,
$\overline{\mbox{Im}\,a}$ is an infinite-dimensional subspace and hence $p$ (as a non-compact operator) cannot belong to $C^*(a)$.
\end{example}

In this light, the following theorem is the best possible extension of Theorem \ref{inv}.

\begin{theorem}\label{invinfty}
Let $a$ be a positive element in a $C^*$-algebra $\A \subseteq \Bbb B({\mathscr H})$ with an infinite spectrum. Then $\lim_{m\rightarrow \infty}ap_m(a)a=a$. In particular, if $a$ is an invertible operator, $\lim_{m \rightarrow \infty}p_m(a)=a^{-1}$.
\end{theorem}
\begin{proof}
Since $\sigma(a)$ is not finite, $f_m(a)\neq 0$ for all $m\in\mathbb{N}$, so every $p_m(\lambda)$ is defined by \eqref{peovi}. Take an arbitrary $m\in\mathbb{N}$. For every $\lambda \in\sigma(a)$ we have
\begin{eqnarray*}
  1-\lambda p_m(\lambda) &=& \left(1-\frac{\lambda}{\|a\|}\right)- \sum_{l=1}^m\frac{\lambda\prod_{k=0}^{l-1}g_k(\lambda)^2}{\prod_{k=0}^l\|f_k(a)\|}\\
   &=& \frac{g_0(\lambda)}{\|a\|}- \sum_{l=1}^m\frac{f_l(\lambda)\prod_{k=0}^{l-1}g_k(\lambda)}{\prod_{k=0}^l\|f_k(a)\|} \\
   &=& \left(\frac{g_0(\lambda)}{\|a\|}-\frac{f_1(\lambda)g_0(\lambda)}{\|a\||f_1(a)\|}\right)- \sum_{l=2}^m\frac{f_l(\lambda)\prod_{k=0}^{l-1}g_k(\lambda)}{\prod_{k=0}^l\|f_k(a)\|} \\
 &=& \frac{g_0(\lambda)g_1(\lambda)}{\|a\|\|f_1(a)\|}- \sum_{l=2}^m\frac{f_l(\lambda)\prod_{k=0}^{l-1}g_k(\lambda)}{\prod_{k=0}^l\|f_k(a)\|}\\
& =&\ldots \\
 &=& \frac{\prod_{k=0}^{m-1}g_k(\lambda)}{\prod_{k=0}^{m-1}\|f_k(a)\|}- \frac{f_m(\lambda)\prod_{k=0}^{m-1}g_k(\lambda)}{\prod_{k=0}^m\|f_k(a)\|}\\
 &=& \frac{\prod_{k=0}^m g_k(\lambda)}{\prod_{k=0}^m\|f_k(a)\|}.
\end{eqnarray*}
Then $$\lambda-\lambda^2 p_m(\lambda)=\frac{\lambda\prod_{k=0}^m g_k(\lambda)}{\prod_{k=0}^m\|f_k(a)\|}=\frac{f_{m+1}(\lambda)}{\prod_{k=0}^m\|f_k(a)\|}$$
and therefore
\begin{eqnarray*}
\|a-ap_m(a)a\|&=&\sup_{\lambda\in\sigma(a)}\{|\lambda-\lambda^2 p_m(\lambda)|\}\\
&=&\sup_{\lambda\in\sigma(a)}\{\frac{\left|f_{m+1}(\lambda)\right|}{\prod_{k=0}^m\|f_k(a)\|}\}\le\frac{\|f_{m+1}(a)\|}{\prod_{k=0}^m\|f_k(a)\|}.
\end{eqnarray*}
From $\sigma(f_k(a))\subseteq[0,\|f_k(a)\|]$ and $f_{k+1}(\lambda)=\|f_{k}(a)\|f_{k}(\lambda)-f_{k}(\lambda)^2$ it follows that $\sigma(f_{k+1}(a))\subseteq[0,\frac{1}{4}\|f_k(a)\|^2],$ so $\|f_{k+1}(a)\|\le \frac{1}{4}\|f_k(a)\|^2$ for all $k.$  Then
\begin{eqnarray*}
\|a-ap_m(a)a\|&\le&\frac{\|f_{m+1}(a)\|}{\prod_{k=0}^m\|f_k(a)\|}\le \frac{1}{4} \frac{\|f_m(a)\|^2}{\prod_{k=0}^m\|f_k(a)\|}=\frac{1}{4} \frac{\|f_m(a)\|}{\prod_{k=0}^{m-1}\|f_k(a)\|}\\
&\le &\frac{1}{4^2} \frac{\|f_{m-1}(a)\|^2}{\prod_{k=0}^{m-1}\|f_k(a)\|}=
\frac{1}{4^2} \frac{\|f_{m-1}(a)\|}{\prod_{k=0}^{m-2}\|f_k(a)\|}\le\ldots\\
&\le& \frac{1}{4^m} \frac{\|f_1(a)\|}{\|a\|}=\frac{1}{4^{m+1}} \frac{\|a\|^2}{\|a\|}=\frac{1}{4^{m+1}}\|a\|.
\end{eqnarray*}
Since $m$ is arbitrary, we conclude that $\lim_{m\to\infty}ap_m(a)a=a.$
\end{proof}

\begin{remark}\label{jako}
From $\lim_{m\rightarrow \infty}ap_m(a)a=a$ one easily gets $\lim_{m\rightarrow \infty}ap_m(a)=p$ in the strong operator topology (where, as before, $p$ denotes the orthogonal projection to $\overline{\mbox{Im}\,a}$).
\end{remark}

By the preceding remark, $p$ is the only possible norm-limit of the sequence $(ap_m(a))$.
In the following proposition we characterize those positive operators $a$ for which the sequence $(ap_m(a))$ converges to $p$ in norm. First we need a lemma. Keeping the notation from the preceding paragraphs, let us also fix the following notational conventions:
for a positive operator $a \in \Bbb B({\mathscr H})$ on a Hilbert space ${\mathscr H}$ denote ${\mathscr H}_1=\overline{\mbox{Im}\,a}$ and ${\mathscr H}_2=\mbox{Ker}\,a$. According to the decomposition ${\mathscr H}={\mathscr H}_1 \oplus {\mathscr H}_2$ we can write $a=\left[\begin{array}{cc}a_1&0\\0&0\end{array}\right]$. For the operators $a$ and $a_1$ we denote by $(f_m(a))$ and $(f_m^{(1)}(a_1))$ the sequences defined by (\ref{efovi}) and by
$(p_m(a))$ and $(p_m^{(1)}(a_1))$ those defined by (\ref{peovi}) and (\ref{peovi2}).

\begin{lemma}
$f_m(a)=\left[\begin{array}{cc}f_m^{(1)}(a_1)&0\\0&0\end{array}\right]$ and
$p_m(a)=\left[\begin{array}{cc}p_m^{(1)}(a_1)&0\\0&0\end{array}\right]$, $\forall m \geq 0$.
\end{lemma}
\begin{proof}
The first assertion is trivial for $m=0$. Let $I_j$ denote the identity operator on ${\mathscr H}_j$ for $j=1,2$. Observe that $\|a\|=\|a_1\|$ which means $\|f_0(a)\|=\|f_0^{(1)}(a_1)\|$. This implies
$f_1(a)=\left[\begin{array}{cc}a_1&0\\0&0\end{array}\right]
\left[\begin{array}{cc}\|f_0(a)\|I_1-f_0(a_1)&0\\0&\|f_0(a)\|I_2\end{array}\right]=
\left[\begin{array}{cc}a_1(\|f_0^{(1)}(a_1)\|I_1-a_1)&0\\0&0\end{array}\right]=
\left[\begin{array}{cc}f_1^{(1)}(a_1)&0\\0&0\end{array}\right]$.
A general inductive argument is obtained exactly in the same way.

The second assertion now follows from the first one combined with (\ref{pmfm}).
\end{proof}

\begin{proposition}\label{closed}
Let $a\in\Bbb{B}({\mathscr H})$ be a positive operator and $p\in\Bbb{B}({\mathscr H})$ the orthogonal projection to
$\overline{\mbox{Im}\,a}.$ Then $(ap_m(a))_m$ converges to $p$ in norm if and only if $\mbox{Im}\,a$
is a closed subspace of ${\mathscr H}$.
\end{proposition}
\begin{proof}
Suppose first that $a$ has a closed range, i.e.~$\mbox{Im}\,a=\overline{\mbox{Im}\,a}$. Then $a_1=a_{|\mbox{Im}\,a}:\mbox{Im}\,a \rightarrow \mbox{Im}\,a$ is a bijection. Since $\mbox{Im}\,a$ is a Hilbert space, $a_1$ is an invertible operator. By Theorem~\ref{invinfty}, the sequence
$(a_1p_m^{(1)}(a_1))$ converges in norm and $\lim_{m \rightarrow \infty}a_1p_m^{(1)}(a_1)=I_1$. By the preceding lemma
$ap_m(a)=\left[\begin{array}{cc}a_1p_m^{(1)}(a_1)&0\\0&0\end{array}\right]$ converges in norm to $\left[\begin{array}{cc}I_1&0\\0&0\end{array}\right]$ which is
the orthogonal projection to $\mbox{Im}\,a=\overline{\mbox{Im}\,a}$.

Conversely, suppose that $(ap_m(a))$ converges in norm. As we already noted, the limit is then necessarily $p$. By the second assertion of the preceding lemma, $I_1$ is then the norm-limit of the sequence $(a_1p_m^{(1)}(a_1))$. Since the group of invertible operators is open, it follows that $a_1p_m^{(1)}(a_1)$ is an invertible operator, for $m$ large enough. In particular, $a_1p_m^{(1)}(a_1)$ is a surjection and hence $a_1$ is a surjection as well. Thus, $\mbox{Im}\,a_1=\overline{\mbox{Im}\,a}$. Since, obviously, $\mbox{Im}\,a=\mbox{Im}\,a_1$, this shows that $a$ has a closed range.
\end{proof}

\vspace{.1in}

Notice that each positive operator with a finite spectrum has a closed range. Thus, the preceding proposition is in the accordance with Theorem~\ref{inv}. At the same time, it provides another explanation of Example \ref{compact} since a compact positive operator with an infinite spectrum cannot have a closed range.

\vspace{.1in}

{\em Concluding remarks:}
(a) To complete our analysis, let us first turn back to the sequence of inequalities from Theorem \ref{induct}.

If $z \in \X$ has the property that $\sigma(\ip{z}{z})$ is finite then, by Proposition \ref{kon} and Theorem \ref{inv}, there exists $M \in \Bbb N$ such that
$f_M(\ip{z}{z})=0$, $f_{M+1}(\ip{z}{z})\not =0$ and $\ip{z}{z}p_M(\ip{z}{z})$ is the projection to $\textup{Im}\ip{z}{z}.$ In this case, the sequence of inequalities from Theorem \ref{induct} is finite and the last term between $\frac{1}{\|z\|^2}\ip{x}{z}\ip{z}{x}$ and $\langle x,x\rangle$ is
$\ip{x}{z} p_M(\ip{z}{z})\ip{z}{x}$ (for all $x$). The following claim explains the reason: the sequence terminates at that place because
$\ip{x}{z} p_M(\ip{z}{z})\ip{z}{x}$ is the maximal element of the set of positive elements under consideration.

{\bf Claim.} {\em Let $\X$ be a semi-inner product module over a $C^*$-algebra $\A\subseteq \Bbb{B}({\mathscr H})$. For $z\in\X$ and $a=\la z,z\ra\in\A$, let $p\in\Bbb{B}({\mathscr H})$ denote the orthogonal projection to
$\overline{\mbox{Im}\,a}.$
Suppose that there exists a positive operator $h\in \Bbb{B}({\mathscr H})$ such that for all $x\in\X$ and every $m\ge 0$ it holds
\begin{equation}\label{ha}
\ip{x}{x} \ge \ip{x}{z} h\ip{z}{x}\ge \ip{x}{z} p_m(\ip{z}{z})
\ip{z}{x}.
\end{equation}
Then $aha=a$ and $ah=p.$}
\begin{proof}
It follows from \eqref{ha} that
$$\la z,z\ra\ge \la z,z\ra h\la z,z\ra\ge\la z,z\ra p_m(\la z,z\ra)\la z,z\ra,\quad \forall m\ge 0,$$
that is, $a\ge aha\ge ap_m(a)a$ for all $m\ge 0.$ By Theorem~\ref{invinfty} (or Remark~\ref{jako}), it follows that $aha=a.$ This implies $ah=p$.
\end{proof}

Suppose now, as in the discussion preceding the above claim, that there exists $M \in \Bbb N$ such that $f_M(a)\not=0$ and $f_{M+1}(a)=0$. Then $a$ has a finite spectrum, $\mbox{Im}\,a$ is a closed subspace, and
$ap_M(a)=p$. So, if $h$ is as in the above claim, then $ah=p$ and therefore $ap_M(a)=ah$. By taking adjoints we get $ha=p_M(a)a$ and this shows that $h$ and $p_M(a)$ coincide on $\mbox{Im}\,a$.

If $\sigma(a)$ is infinite, there is no $M$ as above, but still the sequence $(ap_m(a)a)$ converges in norm to $a$. From the proof of the claim it follows that for any $h\in \Bbb{B}({\mathscr H})$ which satisfies left-hand side inequality of \eqref{ha} it holds $aha\le a.$ Therefore, we have a kind of best result even in this case, since $h$ which appears in \eqref{ha} is such that  $aha=\lim_{m\to\infty}ap_m(a)a=a$.

\vspace{.1in}

(b) Observe that, if $b\in \A$ is positive and such that $\|zb^{\frac{1}{2}}\|\le 1,$ then, by Theorem \ref{p1}, we have
$$\ip{x}{x}\ge \|zb^{\frac{1}{2}}\|^2\ip{x}{x}\ge \ip{x}{zb^{\frac{1}{2}}}\ip{zb^{\frac{1}{2}}}{x}=\ip{x}{z}b\ip{z}{x}$$
for every $x\in \X.$ Thus, the inequalities from Theorem
\ref{induct} can alternatively be derived from the inequalities
$\|zp_m(\ip{z}{z})^{\frac{1}{2}}\|\le 1,\,m \in \Bbb N$. Instead of
proving these inequalities directly, we opted for the inductive
approach from the proof of Theorem \ref{induct} since it leads
naturally to the sequence $(f_m(a))$ and gives us more insight into the
sequence $p_m(a)$ which, as we have seen, has many interesting
properties.

(c) The assertion of Theorem \ref{induct} can be formulated for Gram
matrices as well; the proof requires no essential changes. In this
way, one obtains the result that directly improves Theorem \ref{p1}.\\

\textbf{Acknowledgment.} The authors wish to thank the referee for several helpful comments and valuable suggestions.
The third author was supported by a grant from Ferdowsi University of Mashhad (No. MP87037MOS).


\end{document}